\documentclass[10pt]{amsart}

\usepackage{amsmath}
\usepackage{amssymb}
\usepackage{latexsym}
\usepackage{amsthm}

\newcommand{\FT}{\,\, \widehat{} \,\,}
\newcommand{\stFT}{\,\, \widetilde{} \,\,}

\newcommand{\abs}[1]{\left\vert #1 \right\vert}
\newcommand{\bigabs}[1]{\bigl\vert #1 \bigr\vert}

\newcommand{\norm}[1]{\left\Vert #1 \right\Vert}
\newcommand{\bignorm}[1]{\bigl\Vert #1 \bigr\Vert}

\newcommand{\Sobnorm}[2]{\norm{#1}_{H^{#2}}}

\newcommand{\mixednorm}[3]{\norm{#1}_{L_{t}^{#2}L_{x}^{#3}}}

\newcommand{\C}{\mathbb{C}}

\newcommand{\R}{\mathbb{R}}

\newcommand{\charfn}{\mathbf{1}}

\newcommand{\innerprod}[2]{\left\langle \, #1 , #2 \, \right\rangle}
\newcommand{\biginnerprod}[2]{\bigl\langle \, #1 , #2 \, \bigr\rangle}

\newcommand{\elwt}[1]{\left\langle #1 \right\rangle}

\DeclareMathOperator{\im}{Im}

\newtheorem{theorem}{Theorem}

\newtheorem{lemma}{Lemma}

\theoremstyle{definition}
\newtheorem{definition}{Definition}

\theoremstyle{remark}
\newtheorem{remark}{Remark}

\numberwithin{equation}{section}

\title[Null structure and almost optimal local regularity for DKG]{Null structure and almost optimal local regularity for the Dirac-Klein-Gordon system}

\author[P. D'Ancona]{Piero D'Ancona}

\address{Department of Mathematics\\
University of Rome ``La Sapienza''\\
Piazzale Aldo Moro 2\\
I-00185 Rome\\ Italy}
\email{dancona@mat.uniroma1.it}

\author[D. Foschi]{Damiano Foschi}
\address{Department of Pure and Applied Mathematics\\
University of L'Aquila\\
Via Vetoio, loc. Coppito\\
I-67010 L'Aquila\\ Italy}
\email{foschi@univaq.it}

\author[S. Selberg]{Sigmund Selberg}
\address{Department of Mathematical Sciences\\
Norwegian University of Science and Technology\\
Alfred Getz' vei 1\\
N-7491 Trondheim\\ Norway}
\thanks{The last author was supported by the Research Council of Norway, project no.\ 160192/V30, ``PDE and Harmonic Analysis''.}
\email{sigmund.selberg@math.ntnu.no}


\begin{document}

\begin{abstract}
 We prove almost optimal local well-posedness for the coupled
 Dirac-Klein-Gordon (DKG) system of equations in $1+3$ dimensions.
 The proof relies on the null structure of the system, combined with
 bilinear spacetime estimates of Klainerman-Machedon type.
 It has been known for some time that the Klein-Gordon part of the
 system has a null structure; here we uncover an additional
 null structure in the Dirac equation, which cannot be seen
 directly, but appears after a duality argument.
\end{abstract}

\maketitle

\section{Introduction}

In standard notation, the coupled Dirac-Klein-Gordon (DKG) system
of equations on $\R^{1+3}$ reads
\begin{equation}\label{DKGv1}
\left\{
\begin{aligned}
  &\left( -i\gamma^\mu \partial_\mu + M \right) \psi = g \phi
  \psi,
  \qquad \left( M \ge 0, \, g > 0 \right)
  \\
  &\bigl( -\square + m^2 \bigr) \phi = g \psi^\dagger \gamma^0
  \psi,
  \qquad \bigl( \square = -\partial_t^2 + \Delta, \, m \ge 0 \bigr)
\end{aligned}
\right.
\end{equation}
where the unknowns are (i) a spinor field $\psi(t,x) \in \C^4$,
regarded as a column vector in $\C^4$, and (ii) a real scalar
field $\phi(t,x)$. We use coordinates $t = x^0, \, x =
(x^1,x^2,x^3)$ on $\R^{1+3}$, and write $\partial_\mu =
\frac{\partial}{\partial x^\mu}$. Greek indices $\mu,\nu$ etc.\
range over $0,1,2,3$, Roman indices $j,k$ etc.\ over $1,2,3$, and
repeated indices are summed over these ranges. Thus, $\gamma^\mu
\partial_\mu = \sum_{\mu = 0}^3 \gamma^\mu \partial_\mu$, where
$\{ \gamma^\mu \}_{\mu = 0}^3$ are the $4 \times 4$ Dirac
matrices, given in $2 \times 2$ block form by
$$
  \gamma^{0} = \begin{pmatrix}
    I & 0  \\
    0 & -I
  \end{pmatrix},
  \qquad
  \gamma^{j} = \begin{pmatrix}
    0 & \sigma^{j}  \\
    -\sigma^{j} & 0
  \end{pmatrix},
$$
where
$$
  \sigma^{1} =
  \begin{pmatrix}
    0 & 1  \\
    1 & 0
  \end{pmatrix},
  \qquad
  \sigma^{2} =
  \begin{pmatrix}
    0 & -i  \\
    i & 0
  \end{pmatrix},
  \qquad
  \sigma^{3} =
  \begin{pmatrix}
    1 & 0  \\
    0 & -1
  \end{pmatrix}
$$
are the Pauli matrices. $\psi^\dagger$ denotes the adjoint, i.e.,
the conjugate transpose, hence
$$
  \psi^\dagger \gamma^0 \psi \equiv \abs{\psi_1}^2 +
  \abs{\psi_2}^2 - \abs{\psi_3}^2 - \abs{\psi_4}^2,
$$
where $\psi_1,\dots,\psi_4$ are the components of $\psi$. The
following related matrices occur frequently in the sequel:
$$
  \beta \equiv \gamma^0,
  \qquad
  \alpha^{j} \equiv \gamma^{0}\gamma^{j}
  = \begin{pmatrix}
    0 & \sigma^{j}  \\
    \sigma^{j} & 0
  \end{pmatrix},
  \qquad
  S^{m} \equiv i \gamma^{k} \gamma^{l} =
  \begin{pmatrix}
    \sigma^{m} & 0  \\
    0 & \sigma^{m}
  \end{pmatrix},
$$
where $(k,l,m)$ is any cyclic permutation of $(1,2,3)$. Note the
identities
\begin{align}
  \label{MatrixIdentities1}
  \alpha^j \beta &= - \beta \alpha^j,
  \\
  \label{MatrixIdentities2}
  \alpha^{j} \alpha^{k} &= - \alpha^{k} \alpha^{j} + 2 \delta^{jk} I
  \\
  \label{MatrixIdentities3}
  \alpha^{j} \alpha^{k}
  &= \delta^{jk} I + i \epsilon^{jkl} S^{l}.
\end{align}
Also, $\beta^2 = (\alpha^j)^2 = I$ and $\beta^\dagger = \beta$,
$(\alpha^j)^\dagger = \alpha^j$.

Concerning the Cauchy problem, the most fundamental question is
whether global regularity holds, i.e., given smooth, compactly
supported initial data, does DKG have a smooth solution for all
times $t > 0$? For small data, the answer is yes (see
\cite{Bachelot88, Klainerman85,ChoquetBruhat81}), but for large
data it remains an open question, except in the $1+1$ dimensional
case, see Chadam \cite{Chadam73}. In $1+3$ dimensions, global
regularity is known only for a very special class of (large) data:
Chadam and Glassey \cite{ChadamGlassey74} proved it for data
satisfying the constraints $\psi_1(0,x) = \overline{\psi_4}(0,x)$
and $\psi_2(0,x) = - \overline{\psi_3}(0,x)$, which imply that
$\psi^\dagger \gamma^0 \psi$ vanishes initially, and in fact stays
zero in the evolution; later, Bachelot \cite{Bachelot89} extended
this result to cover also small perturbations around such data.
Another global result is proved in \cite{DiasFigueira91} for data
with special symmetry properties.

In order to make progress on the global regularity question, a
natural strategy is to study local (in time) well-posedness (LWP)
for low regularity data, and then try to exploit the conserved
quantities of the system. This strategy was successfully
implemented for the Maxwell-Klein-Gordon (MKG) and Yang-Mills (YM)
equations by Klainerman and Machedon
\cite{KlainermanMachedon94,KlainermanMachedon95a}, who proved LWP
for data with finite energy and then used the conservation of
energy to push this to a global result, thus recovering, in
particular, the classical result of Eardley and Moncrief
\cite{EardleyMoncrief82}. Compared to MKG and YM, however, DKG has
the unpleasant feature that the conserved energy $\int
e(\phi,\psi) \, dx = \text{const.}$ has a density which is not
positive definite (see \cite{GlasseyStrauss79}):
$$
  e(\phi,\psi) = \im \left( \psi^\dagger \alpha^j \partial_j \psi
  \right) - (M-g\phi) \psi^\dagger \beta \psi - \frac{1}{2} \left(
  (\partial_t \phi)^2 + \abs{\nabla \phi}^2 + m^2 \phi^2 \right).
$$
On the other hand, one does have the conservation of charge:
\begin{equation}\label{ConservationOfCharge}
  \int \abs{\psi(t,x)}^2 \, dx = \text{const.}
\end{equation}
which was a key ingredient in Chadam's proof of global
regularity in the $1+1$ dimensional case \cite{Chadam73} (see also
\cite{Bournaveas00,Fang04}).

We are interested in LWP of the Cauchy problem with data
\begin{equation}\label{Data}
  \psi(0,x) = \psi_0(x), \qquad \phi(0,x) = \phi_0(x), \qquad
  \partial_t \phi(0,x) = \phi_1(x)
\end{equation}
with regularity $(\psi_0,\phi_0,\phi_1) \in H^s \times H^r \times
H^{r-1}$ for minimal $s,r \in \R$. Here $H^s = H^s(\R^3)$ is the
Sobolev space with norm
$$
  \norm{f}_{H^s} = \bignorm{\elwt{\xi}^s \widehat f(\xi)}_{L^2_\xi},
$$
where $\widehat f(\xi)$ denotes the Fourier transform of $f(x)$
and $\elwt{\cdot} = 1 + \abs{\cdot}$. We denote by $\dot
H^s$ the corresponding homogeneous space, with norm
$\norm{f}_{\dot H^s} = \bignorm{\abs{\xi}^s \widehat
f(\xi)}_{L^2_\xi}$.

To get an idea of the minimal regularity required for LWP, one can
apply the usual scaling heuristic. In the massless case $M = m =
0$, DKG is invariant under the rescaling
$$
  \psi(t,x) \longrightarrow \frac{1}{L^{3/2}} \psi\Bigl(\frac{t}{L},\frac{x}{L}\Bigr), \qquad
  \phi(t,x) \longrightarrow \frac{1}{L} \phi\Bigl(\frac{t}{L},\frac{x}{L}\Bigr),
$$
hence the scale invariant data space is (in $1+3$ dimensions)
$$
  (\psi_0,\phi_0,\phi_1) \in L^2 \times \dot H^{1/2} \times
\dot H^{-1/2},
$$
and one does not expect well-posedness below this regularity. The
scaling also suggests that $r = 1/2 + s$ is the natural choice.

On the other hand, DKG is a system of nonlinear wave equations
with quadratic nonlinearities
(as can be seen by squaring the
Dirac equation), and for such equations it is well known (see
\cite{Lindblad96}) that, due to nonlinear effects, one cannot hope
to reach the regularity predicted by scaling unless Klainerman's
null condition is satisfied. The null condition is a condition on
the symbol of the quadratic nonlinearities, which cancels the most
dangerous interactions in a product of free waves.

In fact, by classical methods (energy estimates and Sobolev
embeddings) one can prove LWP for data $(\psi_0,\phi_0,\phi_1) \in
H^{1+\varepsilon} \times H^{3/2+\varepsilon} \times
H^{1/2+\varepsilon}$ for any $\varepsilon > 0$. This can be
improved to $H^{1/2+\varepsilon} \times H^{1+\varepsilon} \times
H^{\varepsilon}$ by using Strichartz type estimates for the
homogeneous wave equation (see \cite{PonceSideris,Bournaveas99}),
but in order to lower the regularity further one needs null
structure. Klainerman and Machedon \cite{KlainermanMachedon92} demonstrated a
null structure, via an auxiliary variable, in the quadratic form
$\psi^\dagger \gamma^0 \psi$ appearing in the Klein-Gordon equation.
Later, Beals and Bezard \cite{BealsBezard96} found a more direct expression of this null
structure by using the eigenspace projections of the Dirac operator
(thus avoiding the auxiliary variable), and they applied it to prove a
smoothing estimate for $\phi$. Specifically, from the
$H^{1/2+\varepsilon} \times H^{1+\varepsilon} \times
H^{\varepsilon}$ result one gets also LWP in $H^{1} \times H^{3/2}
\times H^{1/2}$ by persistence of higher regularity. Then applying
Strichartz type estimates to the equation for $\phi$ in
\eqref{DKGv1} one can see that in fact $\phi(t) \in
H^{2-\varepsilon}$ for every $\varepsilon > 0$. Using the null
structure, and the bilinear spacetime estimates of Klainerman and
Machedon \cite{KlainermanMachedon93}, Beals and Bezard improved
this to $\phi(t) \in H^2$.

On the other hand, Bournaveas \cite{Bournaveas99}, following the idea of Klainerman
and Machedon \cite{KlainermanMachedon92}, found a null
structure in the Dirac part of the system, and used this to get rid of the epsilon in the
$H^{1/2+\varepsilon} \times H^{1+\varepsilon} \times
H^{\varepsilon}$ result, i.e., he proved LWP in the ``energy
class'' $H^{1/2} \times H^1 \times L^2$. While the null structure found by Bournaveas helps
to a certain extent, it has the drawback that it involves squaring the Dirac equation,
which creates serious difficulties at very low regularity.
It should be noted, however, that Bournaveas did not
make use of $X^{s,b}$ type spaces (although they are lurking in
the background in Lemma 2 of \cite{Bournaveas99}), which allow one
to take maximum advantage of the null structure. Using the machinery of these spaces together with the null structure proved in \cite{BealsBezard96, Bournaveas99} and bilinear spacetime estimates of Klainerman-Machedon type, one can, not surprisingly, improve the result from \cite{Bournaveas99}.
In fact, quite recently Fang and Grillakis \cite{FangGrillakis05} have proved LWP in $H^s \times H^1 \times L^2$ for all $1/4 < s \le 1/2$.

The new idea which drives the present paper is that the null form
$\psi^\dagger \gamma^0 \psi$ occurs not only in the Klein-Gordon part,
but in fact also in the Dirac part of the system, as can be seen
via a duality argument. Using this new fact we can get arbitrarily
close to the scale invariant regularity:

\begin{theorem}\label{Thm1}
DKG in $1+3$ dimensions is LWP for data
$$
  (\psi_0,\phi_0,\phi_1) \in H^\varepsilon \times H^{1/2+\varepsilon} \times H^{-1/2+\varepsilon}
$$
for all $\varepsilon > 0$.
\end{theorem}

Note that Theorem \ref{Thm1} leaves open the critical case
($\varepsilon = 0$). One may hope that DKG is globally well posed
for small data in some Besov norm with the same scaling as $L^2
\times \dot H^{1/2} \times \dot H^{-1/2}$, but we do not consider
this question here.

In a forthcoming paper we prove analogous results for the Maxwell-Dirac system, i.e., almost optimal LWP
in $1+3$ dimensions.

\section{Preliminaries}

For convenience we rewrite \eqref{DKGv1} in a slightly different
form, multiplying the Dirac equation on the left by $\beta =
\gamma^0$ to get
\begin{equation}\label{DKGv2}
\left\{
\begin{aligned}
  &-i \left( \partial_t + \alpha \cdot \nabla \right) \psi = -M\beta\psi + \phi
  \beta\psi,
  \\
  &\square\phi = m^2 \phi - \innerprod{\beta\psi}{\psi}_{\C^4},
\end{aligned}
\right.
\end{equation}
where we have also set $g=1$. Here $\alpha$ denotes the vector
$(\alpha^1,\alpha^2,\alpha^3)$ whose components are the Dirac
matrices $\alpha^j=\gamma^0\gamma^j$; thus, $\alpha \cdot \nabla =
\alpha^j \partial_j$. Further, $\innerprod{\cdot}{\cdot}_{\C^4}$ denotes
the standard inner product on $\C^4$.

The operator $-i(\partial_t + \alpha \cdot \nabla)$ is rather
complicated, since $-i\alpha \cdot \nabla$ mixes the components of
the spinor it acts on. To simplify matters, it is natural to
diagonalize by decomposing the spinor field relative to an
eigenbasis of the operator $-i\alpha\cdot\nabla$. The symbol of
the latter is $\alpha \cdot \xi\,$ ($\xi \in \R^3$). A quick
calculation using \eqref{MatrixIdentities2} gives $(\alpha \cdot
\xi)^2 = \abs{\xi}^2 I$, hence the eigenvalues of $\alpha \cdot
\xi$ are $\pm \abs{\xi}$. By symmetry, each eigenspace is
two-dimensional, and the projections onto these eigenspaces are
given by
\begin{equation}\label{ProjectionDef}
  \Pi_\pm(\xi) = \frac{1}{2} \left( I \pm \hat \xi \cdot \alpha
  \right) \qquad \text{where} \qquad \hat \xi \equiv
  \frac{\xi}{\abs{\xi}}.
\end{equation}
Now write
$$
  \psi = \psi_+ + \psi_- \qquad \text{where} \qquad \psi_\pm =
  \Pi_\pm(D) \psi.
$$
Here $D = \nabla/i$, which has Fourier symbol $\xi$. Throughout we use the notation $h(D)$ for the multiplier with symbol $h(\xi)$, for a given function $h : \R^3 \to \C$.
 
Applying $\Pi_\pm(D)$ to the Dirac equation in \eqref{DKGv2}, and
using the identities
$$
  -i\alpha\cdot\nabla = \abs{D}\Pi_+(D) - \abs{D}\Pi_-(D),
$$
and
\begin{equation}\label{ProjectionIdentity1}
  \Pi_\pm(\xi) \beta = \beta \Pi_\mp(\xi)
\end{equation}
(the latter due to \eqref{MatrixIdentities1}), we obtain
\begin{equation}\label{DKGv3}
\left\{
\begin{aligned}
  &\left( -i\partial_t \pm \abs{D} \right) \psi_\pm = -M\beta\psi_\mp + \Pi_\pm(D) \left(\phi
  \beta\psi\right),
  \\
  &\square\phi = m^2 \phi - \innerprod{\beta\psi}{\psi}_{\C^4},
\end{aligned}
\right.
\end{equation}
which is the system we shall work with.

We iterate $\psi_\pm$ and $\phi$ in $X^{s,b}$ type spaces associated to the operators
$-i\partial_t \pm \abs{D}$ and $\square$, whose symbols are $\tau
\pm \abs{\xi}$ and $\tau^2-\abs{\xi}^2$, respectively. The
notation $\widetilde u(\tau,\xi)$ is used for the spacetime Fourier
transform of a function $u(t,x)$.

\begin{definition}
Let $X_\pm^{s,b}$ ($s,b \in \R$) be the completion of the Schwartz
space $\mathcal{S}(\R^{1+3})$ with respect to the norm
$$
  \norm{u}_{X_\pm^{s,b}} = \bignorm{ \elwt{\xi}^s
  \elwt{\tau\pm\abs{\xi}}^b \widetilde
  u(\tau,\xi)}_{L^2_{\tau,\xi}},
$$
where as before $\elwt{\cdot} = 1 + \abs{\cdot}$. Note that
$\norm{u}_{X_\pm^{s,b}} = \bignorm{ \elwt{D}^s
  \elwt{-i\partial_t \pm \abs{D}}^b u}_{L^2_{t,x}}$, by
Plancherel's theorem.
\end{definition}

Spaces of this type were first used by Bourgain \cite{Bourgain93}
for periodic solutions of nonlinear Schr\"odinger and KdV
equations, and later by Kenig, Ponce and Vega
\cite{KenigPonceVega94} in the nonperiodic case. Similar spaces
for the wave equation were first used by Klainerman and Machedon
\cite{KlainermanMachedon95}, who used the notation $H^{s,b}$. Here
we rely on a slight variation of the $H^{s,b}$ spaces of
Klainerman and Machedon, introduced in \cite{Selberg99}
(alternatively, see \cite{Selberg02}) and applied in
\cite{KlainermanSelberg02}, where they are referred to as
wave-Sobolev spaces. To describe these spaces, it is convenient to
introduce mulipliers $D_\pm$ with symbols
$\abs{\tau}\pm\abs{\xi}$. Thus, $\square = D_+ D_-$.

\begin{definition}
  Let $H^{s,b}$ and $\mathcal H^{s,b}$ ($s,b \in \R$) be the
  completions of $\mathcal S(\R^{1+3})$ with respect to the norms
  \begin{align*}
    \norm{u}_{H^{s,b}} &= \bignorm{\elwt{D}^s \elwt{D_-}^b
    u}_{L^2_{t,x}}
    = \bignorm{
    \elwt{\xi}^s \elwt{\abs{\tau} - \abs{\xi}}^b \widetilde
    u(\tau,\xi)}_{L^2_{(\tau,\xi)}},
    \\
    \norm{u}_{\mathcal H^{s,b}} &= \norm{u}_{H^{s,b}} + \norm{\partial_t
    u}_{H^{s-1,b}},
  \end{align*}
  respectively. Observe that the last norm is equivalent to
  $\bignorm{\elwt{D}^{s-1} \elwt{D_+} \elwt{D_-}^b
    u}_{L^2_{t,x}}$.
\end{definition}

We shall also need the restrictions of these spaces to a time slab
$$
  S_T = (0,T) \times \R^3,
$$
since we study local in time solutions. The restriction
$X_\pm^{s,b}(S_T)$ is a Banach space with respect to the norm
$$
  \norm{u}_{X_\pm^{s,b}(S_T)} = \inf \left\{
  \norm{v}_{X_\pm^{s,b}} : \text{$v \in X_\pm^{s,b}$ and $v = u$ on $S_T$}
  \right\}.
$$
In fact, the completeness follows from a basic result of abstract
functional analysis, since $X_\pm^{s,b}(S_T)$ is nothing else than
the quotient space $X_\pm^{s,b}/\mathcal M_\pm$, where $\mathcal
M_\pm$ is the closed subspace $\{ v \in X_\pm^{s,b} : \text{$v =
0$ on $S_T$} \}$. The restriction spaces $H^{s,b}(S_T)$ and
$\mathcal H^{s,b}(S_T)$ are defined analogously.

\section{Null structure}

In this section we discuss the null structure in DKG. First, however,
let us recall the null condition of Klainerman and give a
heuristic argument showing its significance for regularity of
nonlinear waves. To this end, consider a nonlinear wave equation
with a quadratic nonlinearity, $\square u = B(u,u)$, where $B$ is
a bilinear operator given by a Fourier symbol $b$. Specifically,
if $X = (\tau,\xi)$, $Y = (\lambda,\eta)$ and $Z = (\mu,\zeta)$
are vectors in Fourier space $\R \times \R^3$, $B$ is of the form
\begin{equation}\label{Bdef}
\begin{split}
  [B(v,v)]\stFT(X) &= \iint_{Y+Z = X} b(Y,Z) \widetilde v(Y)
  \widetilde v(Z) \, dY \, dZ
  \\
  &= \int b(Y,X-Y) \widetilde v(Y)
  \widetilde v(X-Y) \, dY.
\end{split}
\end{equation}
We say $X = (\tau,\xi)$ is \emph{null} if it lies on the null cone
(light cone) $\abs{\tau} = \abs{\xi}$; this is equivalent to
saying that the symbol $\square(X) \equiv \tau^2-\abs{\xi}^2$ of
the wave operator vanishes on $X$. Let us suppose $v$ is a free
wave, $\square v = 0$, so that $\widetilde v$ is a measure
supported on the null cone, and let us look at the regularity of
$u$ solving $\square u = B(v,v)$. (This problem arises naturally
when solving the nonlinear problem by iteration.) In Fourier
space,
$$
  \square(X) \widetilde u(X) = [B(v,v)]\stFT(X),
$$
so one gains a lot of regularity when $X$ is away from the null
cone ($\abs{\square(X)} \gtrsim 1$). Near the cone, things are not
so favorable, but if it happens that $[B(v,v)]\stFT(X)$ vanishes (to
some order) when $X$ is null, this should improve the regularity
in this difficult region. But $X = Y+Z$, where $Y,Z$ are null
(since now $v$ in \eqref{Bdef} is a free wave), hence $X$ is
null if and only if $Y,Z$ are parallel. One concludes that
$[B(v,v)]\stFT(X)$ vanishes for null $X$ if Klainerman's null
condition is satisfied, i.e., if
\begin{equation}\label{NullCondition}
  b(Y,Z) = 0 \quad \text{for} \quad Y,Z \quad \text{null
  parallel.}
\end{equation}

\begin{remark}\label{Remark1}
Note that if two null vectors $Y=(\lambda,\eta)$ and
$Z=(\mu,\zeta)$ are on the same component of the cone, i.e., if
$\lambda$ and $\mu$ have the same sign, then they are parallel if
and only if $\angle(\eta,\zeta)=0$, whereas if they are on
opposite components of the cone, the condition is
$\angle(\eta,-\zeta) = 0$. Here and in the sequel we use the notation $\angle(\eta,\zeta)$
for the angle between two vectors in $\R^3$.
\end{remark}

Let us now turn to the null structure in DKG, starting with the
Klein-Gordon part of the system:
\begin{equation}\label{KG}
  \square \phi = - \innerprod{\beta \psi}{\psi}_{\C^4}.
\end{equation}
For simplicity we set $M = m = 0$ in this section.

Since we deal with spinors, the formulation of the null condition
differs somewhat from the above. First observe that
\begin{equation}\label{MainBilinear}
  [\innerprod{\beta\psi}{\psi}_{\C^4}]\stFT(X) = \iint_{Y+Z = X}
  \biginnerprod{\beta \widetilde\psi(Y)}{\widetilde{\psi}(-Z)}_{\C^4} \, dY \,
  dZ,
\end{equation}
where the minus sign in front of $Z$ stems from the complex
conjugation in the inner product. Second, \eqref{NullCondition}
was derived from the action of \eqref{Bdef} on a free wave
$\square v = 0$, whereas in our case there are two separate
species of free waves, namely $\psi_\pm$ satisfying $(-i\partial_t
\pm \abs{D})\psi_\pm = 0$ (cf.\ \eqref{DKGv3}). Taking data
$\psi(0,x) = \psi_0(x)$ we have
\begin{equation}\label{FreeWaves}
  \psi_\pm(t) = e^{\pm it\abs{D}} \psi_0^\pm \qquad \text{where}
  \qquad \psi_0^\pm = \Pi_\pm(D) \psi_0.
\end{equation}
The spacetime Fourier transforms
\begin{equation}\label{FreeWavesFT}
  \widetilde{\psi_\pm}(Y) = \delta(\lambda\pm\abs{\eta}) \widehat{\psi_0^\pm}(\eta)
  \qquad \left( Y = (\lambda,\eta) \right)
\end{equation}
are supported on opposite components $-\lambda=\pm\abs{\eta}$ of
the null cone $\abs{\lambda} = \abs{\eta}$. Using this information
we can state the null condition for $\innerprod{\beta \psi}{\psi}_{\C^4}$ with $\psi$
replaced by $\psi_\pm$, for all possible combinations of signs (as
before, $Y = (\lambda,\eta)$ and $Z = (\mu,\zeta)$):
\begin{itemize}
  \item[(N1)]
  In the $++$ and $--$ cases, i.e., taking $\innerprod{\beta
  \psi_+}{\psi_+}_{\C^4}$ or $\innerprod{\beta
  \psi_-}{\psi_-}_{\C^4}$ in \eqref{MainBilinear}, we see that $Y, Z$ are on
  opposite components of the cone (because $Y,-Z$ evidently are on the same component), hence the
  null condition (cf.\ Remark \ref{Remark1}) says that the
  (matrix valued) symbol should vanish when the angle $\angle(\eta,-\zeta) = 0$.
  \item[(N2)]
  In the $+-$ and $-+$ cases, $Y, Z$ are on the same component of
  the cone, hence the null condition says that the symbol should
  vanish when the angle $\angle(\eta,\zeta) = 0$.
\end{itemize}

This null condition is indeed satisfied by the symbol of
\begin{equation}\label{MainBilinearSplit}
  (\psi,\psi') \longmapsto
  \innerprod{\beta\Pi_+(D)\psi}{\Pi_\pm(D)\psi'}_{\C^4}.
\end{equation}
This was proved already in \cite{KlainermanMachedon92, BealsBezard96}, but we give a considerably simpler
proof below. The main new contribution in the present paper, however, is the fact that
the null bilinear forms \eqref{MainBilinearSplit} occur not only in the Klein-Gordon part
of the system, but also in the Dirac part. To see this requires a duality argument which we now outline.
To show the main idea unobscured by technical issues, we prefer to present first a heuristic argument
corresponding to the critical regularity $\psi_\pm \in X_\pm^{0,1/2}$;
the rigorous proof of Theorem \ref{Thm1} is then given in the following sections.

In the following heuristic we take zero initial data for $\phi$, so that
$\phi = - \square^{-1} \innerprod{\beta\psi}{\psi}$, where $\square^{-1}F$
denotes the solution of $\square u =F$ with vanishing initial
data. The Dirac part of the system \eqref{DKGv3} then reads (from now on we drop the
index $\C^4$ on the inner product)
$$
  \left( -i\partial_t \pm \abs{D} \right) \psi_\pm = \Pi_\pm(D) F
  \qquad \text{where}
  \qquad
  F =
  \left( -\square^{-1}\innerprod{\beta\psi}{\psi}\right)
  \beta\psi,
$$
and estimating $\psi_\pm$ in $X_\pm^{0,1/2}$ reduces,
heuristically, to estimating $\norm{\Pi_\pm(D)
F}_{X_\pm^{0,-1/2}}$ (cf.\ Lemma \ref{Lemma3} in the next section). Estimating the
latter by duality, we are led to consider integrals, for spinor
valued $\psi' \in X_\pm^{0,1/2}$,
\begin{align*}
  \iint \innerprod{\Pi_\pm(D)F}{\psi'} \, dt \, dx
  &= \iint \innerprod{F}{\Pi_\pm(D) \psi'} \, dt \, dx
  \\
  &= - \iint \left( \square^{-1}\innerprod{\beta\psi}{\psi}\right)
  \innerprod{\beta\psi}{\Pi_\pm(D) \psi'} \, dt \, dx,
\end{align*}
so indeed, \eqref{MainBilinearSplit} crops up one more time.

In fact, the complete null structure of DKG can be elegantly summed up
in a single line: In the last integral, dropping the prime on $\psi'$ and splitting the other fields using $\Pi_\pm(D)$, we end up with
$$
  \iint \left( \square^{-1}\innerprod{\beta\Pi_\pm(D)\psi}{\Pi_\pm(D)\psi}\right)
  \cdot
  \innerprod{\beta\Pi_\pm(D)\psi}{\Pi_\pm(D)\psi} \, dt \,
  dx,
$$
for all possible combinations of signs. Replacing
$\square^{-1}$ by $\abs{\square}^{-1}$ and distributing it equally
over the two factors (this particular heuristic is based on
Plancherel's theorem) yields
$$
  \iint \left( \abs{\square}^{-1/2}\innerprod{\beta\Pi_\pm(D)\psi}{\Pi_\pm(D)\psi}\right)
  \cdot
  \left( \abs{\square}^{-1/2}\innerprod{\beta\Pi_\pm(D)\psi}{\Pi_\pm(D)\psi}\right) \, dt \,
  dx.
$$
The last integral embodies the complete null structure in DKG, and
shows the striking symmetry of the system. It suggests that the
key problem is to prove the bilinear ``estimate''
\begin{equation}\label{KeyBilinearProblem}
  \norm{\abs{\square}^{-1/2}\innerprod{\beta\Pi_+(D)\psi}{\Pi_\pm(D)\psi}}_{L^2_{t,x}}
  \lesssim
  \norm{\psi}_{X_+^{0,1/2}} \norm{\psi}_{X_\pm^{0,1/2}},
\end{equation}
which fails, but not by much. In fact, we shall reduce Theorem \ref{Thm1}
to certain perturbations around this estimate (see \eqref{4} and \eqref{6}), which in turn reduce, on account
of the null structure, to some bilinear spacetime estimates of the type
first studied by Klainerman and Machedon \cite{KlainermanMachedon93}.
For the free wave case, see \eqref{KeyBilinearProblemFree} below.

Let us now verify that that the null condition (N1), (N2) is satisfied by \eqref{MainBilinearSplit}.
In fact, since $\Pi_\pm(D)$ does not involve time at all, it
suffices to consider spinor fields $\psi(x),\psi'(x)$. Replacing
$X,Y,Z$ in \eqref{MainBilinear} by $\xi,\eta,\zeta \in \R^3$, we
then have
$$
  \left[\innerprod{\beta\Pi_+(D)\psi}{\Pi_\pm(D)\psi'}\right]\FT(\xi)
  = \iint_{\eta+\zeta = \xi}
  \biginnerprod{\beta\Pi_+(\eta)\widehat\psi(\eta)}{\Pi_\pm(-\zeta)\widehat{\psi'}(-\zeta)}
  \, d\eta \, d\zeta,
$$
and since $\Pi_\pm(\xi)^\dagger = \Pi_\pm(\xi)$, we obtain
\begin{equation}\label{SymbolCalculation}
\begin{aligned}
  \biginnerprod{\beta\Pi_+(\eta)\widehat\psi(\eta)}{\Pi_\pm(-\zeta)\widehat{\psi'}(-\zeta)}
  &=
  \biginnerprod{\Pi_\pm(-\zeta)\beta\Pi_+(\eta)\widehat\psi(\eta)}{\widehat{\psi'}(-\zeta)}
  \\
  &=
  \biginnerprod{\beta\Pi_\mp(-\zeta)\Pi_+(\eta)\widehat\psi(\eta)}{\widehat{\psi'}(-\zeta)},
\end{aligned}
\end{equation}
where in the last step we used the commutation identity
\eqref{ProjectionIdentity1}, which inverts the sign.

Thus, we have:

\begin{lemma}\label{SymbolLemma}
The symbol of (\ref{MainBilinearSplit})
is the matrix $\beta \Pi_{\mp}(-\zeta)\Pi_+(\eta)$.
\end{lemma}

The symbol $\beta\Pi_\mp(-\zeta)\Pi_+(\eta)$ does indeed satisfy the null
condition (N1), (N2), by orthogonality of the eigenspaces. In
fact, the symbol vanishes to first order in the angle (note that
the following lemma can be applied in all cases, i.e., for all
combinations of signs, because $\Pi_+(-\xi) = \Pi_-(\xi)$).

\begin{lemma}\label{Lemma5}
$\Pi_+(\xi) \Pi_-(\eta) = O(\theta)$, where $\theta =
\angle(\xi,\eta)$.
\end{lemma}

\begin{proof}
\begin{align*}
  4\Pi_+(\xi) \Pi_-(\eta)
  &= (I+\hat\xi_j\alpha^j)(I-\hat\eta_k\alpha^k)
  \\
  &= I - \hat\xi_j\hat\eta_k\alpha^j\alpha^k +
  (\hat\xi-\hat\eta)\cdot\alpha
  \\
  &= (1-\hat\xi\cdot\hat\eta) I - i\epsilon^{jkl}
  \hat\xi_j\hat\eta_k S^l +
  (\hat\xi-\hat\eta)\cdot\alpha \qquad \left( \text{by
  \eqref{MatrixIdentities3}} \right)
  \\
  &= \underbrace{(1-\hat\xi\cdot\hat\eta)}_{O(\theta^2)} I
  - i \underbrace{(\hat\xi \times \hat\eta)}_{O(\theta)}
  \cdot S +
  \underbrace{(\hat\xi-\hat\eta)}_{O(\theta)}\cdot\alpha
\end{align*}
where $\hat\xi \equiv \xi/\abs{\xi}$ and $S = (S^1,S^2,S^3)$.
\end{proof}

\begin{remark}
For readers familiar with the standard null forms $Q_0$, $Q_{ij}$
and $Q_{0j}$, we point out that the factors
$1-\hat\xi\cdot\hat\eta$, $\hat\xi\times\hat\eta$ and
$\hat\xi-\hat\eta$ are the symbols of
$Q_0(\abs{D}^{-1}u,\abs{D}^{-1}v)$,
$Q_{ij}(\abs{D}^{-1}u,\abs{D}^{-1}v)$ and
$Q_{0j}(\abs{D}^{-1}u,\abs{D}^{-1}v)$, respectively, in the case
of free waves $\square u = \square v = 0$.
\end{remark}

This concludes the discussion of the null structure in DKG.
To illustrate how it is used, let us estimate the left hand side of
\eqref{KeyBilinearProblem} in the important
case of free waves $\psi_\pm$ given by \eqref{FreeWaves},
\eqref{FreeWavesFT}. Taking the $++$ case for the sake of definiteness, we shall prove
\begin{equation}\label{KeyBilinearProblemFree}
  \norm{\abs{\square}^{-1/2}\innerprod{\beta\psi_+}{\psi_+}}_{L^2_{t,x}}
  \le C
  \norm{\psi_0}_{L^2}^2.
\end{equation}
Applying Lemmas \ref{SymbolLemma} and \ref{Lemma5}, we see that
\begin{align*}
  &\abs{\bigl[\abs{\square}^{-1/2}\innerprod{\beta\psi_+}{\psi_+}\bigr]\stFT(\tau,\xi)}
  \\
  &= \frac{1}{\sqrt{\abs{\xi}^2-\tau^2}} \abs{ \int_{\R^3} \delta\bigl( \tau + \abs{\eta} -
  \abs{\eta-\xi} \bigr)
  \innerprod{\beta\Pi_-(\eta-\xi)\Pi_+(\eta)\widehat{\psi_0}(\eta)}{\widehat{\psi_0}(\eta-\xi)}
  \, d\eta }
  \\
  &\le
  \frac{C}{\sqrt{\abs{\xi}^2-\tau^2}} \int_{\R^3} \theta \delta\bigl( \tau + \abs{\eta} -
  \abs{\eta-\xi} \bigr)
  \bigabs{\widehat{\psi_0}(\eta)}\bigabs{\widehat{\psi_0}(\eta-\xi)}
  \, d\eta
  \\
  &\le
  C \int_{\R^3} \delta\bigl( \tau + \abs{\eta} -
  \abs{\eta-\xi} \bigr)
  \frac{\bigabs{\widehat{\psi_0}(\eta)}\bigabs{\widehat{\psi_0}(\eta-\xi)}}
  {\abs{\eta}^{1/2}\abs{\eta-\xi}^{1/2}}
  \, d\eta,
\end{align*}
where $\theta = \angle(\eta,\eta-\xi)$ and we used
$$
  \abs{\xi}^2-\tau^2 = \abs{\xi}^2 - \left( \abs{\eta} -
  \abs{\eta-\xi} \right)^2 = 2\abs{\eta}\abs{\eta-\xi} \left( 1 -
  \cos \theta \right) \approx \abs{\eta}\abs{\eta-\xi} \theta^2.
$$
Here and in the sequel, the notation $X \approx Y$ stands for
$C^{-1} X \le Y \le CX$ where $C>0$ is some absolute constant.

We conclude (going back to physical space and applying H\"older's inequality) that \eqref{KeyBilinearProblemFree} reduces to 
the classical estimate
\begin{equation}\label{StrichartzL4}
  \bignorm{e^{\pm it\abs{D}} f}_{L^4(\R^{1+3})} \le C \norm{f}_{\dot
  H^{1/2}}
\end{equation}
due to Strichartz \cite{Strichartz77}.

\section{Some properties of $X^{s,b}$ and $H^{s,b}$}

Here we recall some basic, well known properties of $X^{s,b}$ and
$H^{s,b}$ spaces, needed in the proof of our main result, Theorem \ref{Thm1}.
For the convenience of the reader we include short
sketches of the proofs in some cases. For more details and further
references, see e.g. \cite{Tao01, KlainermanSelberg02}.

We start with $X^{s,b}$, commenting on the more complicated
$H^{s,b}$ at the end. The discussion applies to $X^{s,b}$ in
general form: Starting from a PDE on $\R^{1+n}$, any $n \ge 1$, of
the form
\begin{equation}\label{DispersivePDE}
  -i\partial_t u = h(D) u
\end{equation}
where $h : \R^n \to \R$ and $h(D)$ is the multiplier with symbol
$h(\xi)$, one defines $X^{s,b}$ ($s,b \in \R$) via the norm
$$
  \norm{u}_{X^{s,b}} = \bignorm{ \elwt{\xi}^s  \elwt{\tau-h(\xi)}^b
  \widetilde u(\tau,\xi)}_{L^2_{\tau,\xi}}.
$$
The cases of interest for us here are $h(\xi) = -\abs{\xi}$, which
gives $X_+^{s,b}$, and $h(\xi) = \abs{\xi}$, which gives
$X_-^{s,b}$, but we prefer to keep the general notation (and
general dimension) in this section.

Note that $e^{ith(D)}$ is the free propagator of
\eqref{DispersivePDE}. The first important observation is that the
elements of $X^{s,b}$ are superpositions of free solutions with
$H^s$ data, suitably modulated:

\begin{lemma}\label{Lemma1}
$u \in X^{s,b}$ if and only if there exists $f \in L^2 \bigl(
\elwt{\lambda}^{2b} \, d\lambda; H^s(\R^n) \bigr)$ such that
\begin{equation}\label{Superposition}
  u(t) = \int_\R e^{it\lambda} e^{ith(D)} f(\lambda) \, d\lambda
  \qquad (\text{$H^s$-valued}).
\end{equation}
Moreover, $\norm{u}_{X^{s,b}}^2 = \int_\R
\norm{f(\lambda)}_{H^s}^2 \elwt{\lambda}^{2b} \, d\lambda$.
\end{lemma}

\begin{proof}
The idea is to foliate Fourier space $(\tau,\xi) \in \R \times
\R^n$ by the surfaces $\tau-h(\xi) = \text{const}$. Define
$f(\lambda)$ (a.e.) by $[f(\lambda)]\FT(\xi) = \widetilde
u\bigl(\lambda+h(\xi),\xi\bigr)$. Then
$$
  \widetilde u(\tau,\xi) = \int
  \delta\bigl(\lambda-\tau+h(\xi)\bigr)\widehat{f(\lambda)}(\xi)
  \, d\lambda,
$$
which agrees with the spacetime Fourier transform of
\eqref{Superposition}.
\end{proof}

From this lemma and the dominated convergence theorem for
$H^s$-valued integrals, one easily obtains
\begin{equation}\label{BasicEmbedding}
  X^{s,b} \hookrightarrow C_{\rm b}(\R;H^s) \qquad \text{for} \qquad b
  > \frac{1}{2},
\end{equation}
where $\hookrightarrow$ means continuous inclusion. This in turn
implies $X^{s,b}(S_T) \hookrightarrow C([0,T];H^s)$ for $b > 1/2$.

Another easy but exceedingly useful consequence of Lemma
\ref{Lemma1} is that Strichartz type estimates (linear or
multilinear) for \eqref{DispersivePDE} imply corresponding
estimates for $X^{s,b}$:

\begin{lemma}\label{Lemma2}
Let $T$ be a multilinear operator $(f_{1}(x),\dots,f_{k}(x))
\mapsto T(f_{1},\dots,f_{k})(x)$ acting in $x$-space. If $T$
satisfies an estimate of the form
$$
  \mixednorm{T\bigl(e^{ith_1(D)} f_{1},\dots, e^{ith_k(D)} f_{k}\bigr)}{q}{r}
  \le C \Sobnorm{f_{1}}{s_{1}} \cdots \Sobnorm{f_{k}}{s_{k}},
$$
then
$$
  \mixednorm{T\left(u_{1},\dots, u_{k} \right)}{q}{r}
  \le C_{b}
  \norm{u_{1}}_{X_1^{s_{1},b}} \cdots \norm{u_{k}}_{X_k^{s_{k},b}}
$$
holds for all $u_{j} \in X_j^{s_{j},b}$, $1 \le j \le k$, provided
$b > 1/2$. Here $X_j^{s,b}$ is defined using the symbol
$h_j(\xi)$.
\end{lemma}

\begin{proof}
Since $T$ acts only in $x$, not in $t$, the multilinearity gives,
using the representation from Lemma \ref{Lemma1},
$$
  T\left(u_{1},\dots, u_{k} \right)
  = \int_\R \cdots \int_\R
  e^{it(\lambda_1+\cdots+\lambda_k)} T\bigl(e^{ith_1(D)} f_{1},\dots, e^{ith_k(D)} f_{k}\bigr)
  \, d\lambda_1 \cdots d\lambda_k.
$$
Minkowski's integral inequality followed by Cauchy-Schwarz in
$d\lambda_j$ easily leads to the desired estimate.
\end{proof}

For example, the classical Strichartz estimate \eqref{StrichartzL4}
implies $X_\pm^{1/2,b} \hookrightarrow
L^4(\R^{1+3})$ for $b > 1/2$.

Finally, we consider the LWP of the linear Cauchy problem
associated to \eqref{DispersivePDE},
\begin{equation}\label{AbstractCP}
  -\bigl(i\partial_t + h(D)\bigr)u = F(t,x), \qquad u(0,x) = f(x),
\end{equation}
in the restricted space $X^{s,b}(S_T)$, $b > 1/2$. In the rest of
the paper we impose the condition $T \le 1$ to avoid having to
keep track of the growth of certain constants as $T$ becomes
large.

\begin{lemma}\label{Lemma3}
Let $1/2 < b \le 1$, $s \in \R$, $0 < T \le 1$. Also, let $0 \le
\delta \le 1-b$. Then for all data $F \in X^{s,b-1+\delta}(S_T)$,
$f \in H^s$, the Cauchy problem (\ref{AbstractCP}) has a unique
solution $u \in X^{s,b}(S_T)$, satisfying the first member of
(\ref{AbstractCP}) in the sense of $\mathcal D'(S_T)$. Moreover,
\begin{equation}\label{Xineq}
  \norm{u}_{X^{s,b}(S_T)} \le C \left( \norm{f}_{H^s}
  + T^\delta \norm{F}_{X^{s,b-1+\delta}(S_T)} \right),
\end{equation}
where $C$ only depends on $b$.
\end{lemma}

Note that in Fourier space, $[\tau-h(\xi)] \widetilde u(\tau,\xi)
= \widetilde F(\tau,\xi)$, so heuristically, \eqref{Xineq} says
that in the time localized case, we can replace the singular
symbol $\tau-h(\xi)$ by $\elwt{\tau-h(\xi)}$, and simply divide
out.

In the following proof sketch, we follow closely the argument
given in \cite{KenigPonceVega94}, but with a slight modification
to get the factor $T^\delta$ in the right hand side of
\eqref{Xineq}, which of course is useful in a contraction
argument.

\begin{proof}
Start by picking any extension of $F$, which we still denote $F$.
Then the problem is to prove \eqref{Xineq} without the time
restriction in the right hand side. Split $u = u_0 + u_1$, where
$u_0$ is the homogeneous part, $u_1$ the inhomogeneous part.
Further split $u_1 = u_{1,\text{near}} + u_{1,\text{far}}$
corresponding to the Fourier domains $\abs{\tau-h(\xi)} \le
T^{-1}$ and $\abs{\tau-h(\xi)} > T^{-1}$.

First, $\widetilde{u_0}(\tau,\xi) = \delta\bigl(\tau-h(\xi)\bigr)
\widehat f(\xi)$, hence if $\chi(t)$ is a smooth cut-off such that
$\chi(t) = 1$ for $\abs{t} \le 1$ and $\chi(t) = 0$ for $\abs{t}
\ge 2$, then $\norm{u_0}_{X^{s,b}(S_T)}$ is bounded by
$$
  \norm{\chi(t)
  u_0(t,\cdot)}_{X^{s,b}}
  = \bignorm{ \elwt{\xi}^s \elwt{ \tau-h(\xi)}^b \widehat
  \chi\bigl(\tau-h(\xi)\bigr) \widehat f(\xi)}_{L^2_{\tau,\xi}}
  = \norm{\chi}_{H^b} \norm{f}_{H^s}.
$$

Next, since (we use $\charfn$ to denote the indicator function of
the set determined by the condition in the subscript)
$$
[u_{1,\text{far}}]\stFT(\tau,\xi) =
\frac{\charfn_{\abs{\tau-h(\xi)} > T^{-1}}}{\tau-h(\xi)}
\widetilde F(\tau,\xi),
$$
it follows that $\norm{u_{1,\text{far}}}_{X^{s,b}} \le C T^\delta
\norm{F}_{X^{s,b-1+\delta}}$.

Finally, Duhamel's principle gives $u_{1,\text{near}}(t) =
\sum_{k=1}^\infty \frac{t^k}{k!} e^{ith(D)} g_k,$
where
$\widehat{g_k}(\xi) = \int \left[ i\bigl(\lambda-h(\xi)\bigr)
\right]^{k-1} \charfn_{\abs{\lambda-h(\xi)} \le T^{-1}} \widetilde
F(\lambda,\xi) \, d\lambda$. Using Cauchy-Schwarz in $d\lambda$
one easily verifies that $
  \norm{g_k}_{H^s} \le C T^{b-1/2+\delta-k}
  \norm{F}_{X^{s,b-1+\delta}},
$
so finally,
\begin{align*}
  \norm{u_{1,\text{near}}}_{X^{s,b}(S_T)}
  &\le
  \norm{\chi\left(\frac{t}{T}\right)
  u_{1,\text{near}}(t,\cdot)}_{X^{s,b}}
  \\
  &\le C\sum_{k=1}^\infty \frac{T^k}{k!} \norm{ \left( \frac{t}{T}
  \right)^k \chi\left(\frac{t}{T}\right)}_{H^b_t} \norm{g_k}_{H^s}
  \\
  &\le C\sum_{k=1}^\infty \frac{T^k}{k!} \left( T^{1/2-b}
  \norm{t^k \chi(t)}_{H^b_t} \right) \left( T^{b-1/2+\delta-k}
  \norm{F}_{X^{s,b-1+\delta}} \right)
  \\
  &\le C T^{\delta} \left( \sum_{k=1}^\infty \frac{k2^{k-1}}{k!}
  \right) \norm{F}_{X^{s,b-1+\delta}},
\end{align*}
since $\norm{t^k \chi(t)}_{H^b_t} \le \norm{t^k \chi(t)}_{H^1_t}
\le C_{\chi} \left( 2^k + k2^{k-1} \right)$ by the support
assumption.
\end{proof}

To end this section, let us state the analogues of Lemmas
\ref{Lemma1}--\ref{Lemma3} for $H^{s,b}$. For $u \in H^{s,b}$ we
split $u = u_+ + u_-$ corresponding to the Fourier domains $-\tau
> 0$ and $-\tau < 0$, i.e., we set $\widetilde{u_\pm}(\tau,\xi) =
\charfn_{\pm(-\tau) > 0} \widetilde u(\tau,\xi)$. (Letting $-\tau$
determine the sign is consistent with the choice of signs in the
projections \eqref{ProjectionDef}, since $-\tau$ corresponds to
the energy.) Then $u \in H^{s,b}$ is equivalent to saying that
$u_\pm \in X_\pm^{s,b}$, so Lemma \ref{Lemma1} applies to give a
characterization of $H^{s,b}$ in terms of superimposed free waves.
Moreover, $\norm{u}_{H^{s,b}}^2 = \norm{u_+}_{X_+^{s,b}}^2
+\norm{u_-}_{X_-^{s,b}}^2$, so Lemma \ref{Lemma2} also applies: If
$T$ is as in Lemma \ref{Lemma2} and the estimate
$$
  \mixednorm{T\bigl(e^{\pm it\abs{D}} f_{1},\dots, e^{\pm it\abs{D}} f_{k}\bigr)}{q}{r}
  \le C \Sobnorm{f_{1}}{s_{1}} \cdots \Sobnorm{f_{k}}{s_{k}}
$$
holds, then
$$
  \mixednorm{T\left(u_{1},\dots, u_{k} \right)}{q}{r}
  \le C_{b}
  \norm{u_{1}}_{H^{s_{1},b}} \cdots \norm{u_{k}}_{H^{s_{k},b}}
$$
holds for all $u_{j} \in H^{s_{j},b}$, $1 \le j \le k$, provided
$b
> 1/2$.

Next, consider the Cauchy problem
\begin{equation}\label{WaveCP}
  \square u = F(t,x), \qquad u(0,x) = f(x), \qquad \partial_t
  u(0,x) = g(x).
\end{equation}
One may think that by rewriting \eqref{WaveCP} as a first order
system, Lemma \ref{Lemma3} could be applied directly, but if one
works with data in inhomogeneous Sobolev spaces, there is a
problem with low frequencies. Thus, the following lemma requires a
separate proof, which can be found in \cite{Selberg99, Selberg02};
see \cite{KlainermanMachedon95} for an earlier version of this
estimate (in a slightly different norm, and with $T = 1$).

\begin{lemma}\label{Lemma4}
Let $1/2 < b \le 1$, $s \in \R$, $0 < T \le 1$ and $0 \le \delta
\le 1-b$. Then for all data $F \in H^{s-1,b-1+\delta}(S_T)$, $f
\in H^s$ and $g \in H^{s-1}$, there exists a unique $u \in
\mathcal H^{s,b}(S_T)$ solving (\ref{WaveCP}) on $S_T$. Moreover,
$$
  \norm{u}_{\mathcal H^{s,b}(S_T)} \le C \left( \norm{f}_{H^s} +
  \norm{g}_{H^{s-1}} + \sigma(T)
  \norm{F}_{H^{s-1,b-1+\delta}(S_T)} \right),
$$
where $\sigma > 0$ depends continuously on $T > 0$ and satisfies
$\lim_{T \to 0^+} \sigma(T) = 0$ if $\delta > 0$.
\end{lemma}

\section{Proof of Theorem \ref{Thm1}}

We use the iteration scheme,\footnote{It may be more natural to include the mass terms
in the operators on the left hand sides, but this would require generalizing the Klainerman-Machedon
type estimates to the massive operators; while this can certainly be done (and in some cases has been done),
it is not something we wish to undertake in the present paper. In any event, putting the linear mass terms in the right hand
sides is harmless as far as the local in time contraction is concerned.} for $k \ge -1$,
$$
\left\{
\begin{aligned}
  &\left( -i\partial_t \pm \abs{D} \right) \psi_\pm^{k+1} =
  -M\beta\psi_\mp^k
  + \Pi_\pm(D) \bigl(\phi^k
  \beta\psi^k\bigr),
  \\
  &\square\phi^{k+1} = m^2 \phi^k - \biginnerprod{\beta\psi^k}{\psi^k},
\end{aligned}
\right.
$$
where $\psi^k = \psi_+^k + \psi_-^k$ and $\psi_\pm^{-1}, \phi^{-1}
\equiv 0$. The initial data are
$$
  \psi^k_\pm(0,x) = \Pi_\pm(D)\psi_0(x), \qquad \phi^k(0,x) = \phi_0(x), \qquad
  \partial_t \phi^k(0,x) = \phi_1(x)
$$
where $\psi_0 \in H^s$, $(\phi_0,\phi_1) \in H^{1/2+\varepsilon}
\times H^{-1/2+\varepsilon}$ for $0 < \varepsilon \le 1/2$.
Consider $\varepsilon$ fixed, and let $\varepsilon'>0$ denote a
sufficiently small number, depending on $\varepsilon$. We iterate
in the spaces
$$
  \psi^k_\pm \in X_\pm^{\varepsilon,1/2+\varepsilon'}(S_T),
  \qquad \phi^k \in \mathcal
  H^{1/2+\varepsilon,1/2+\varepsilon'}(S_T),
$$
where $S_T = (0,T) \times \R^3$ for some $0 < T \le 1$ to be
chosen sufficiently small depending on the size of
$$
  \mathcal I_0 = \norm{\psi_0}_{H^\varepsilon}
  + \norm{\phi_0}_{H^{1/2+\varepsilon}}
  + \norm{\phi_1}_{H^{-1/2+\varepsilon}}.
$$
In the estimates that follow, $\lesssim$ stands for $\le$ up to a
multiplicative constant $C$ which may depend on $\varepsilon$ but
is independent of $T$, and $\sigma(T)$ denotes a positive,
continuous function of $0 < T \le 1$ such that $\lim_{T \to 0^+}
\sigma(T) = 0$.

By induction, $\Pi_\pm(D)\psi^k_\pm = \psi^k_\pm$ for all $k$, hence
\begin{equation}\label{0}
  \psi^k = \Pi_+(D)\psi^k_+ + \Pi_-(D)\psi^k_-.
\end{equation}

By Lemmas \ref{Lemma3} and \ref{Lemma4},
\begin{multline}\label{1}
  \bignorm{\psi^{k+1}_\pm}_{X_\pm^{\varepsilon,1/2+\varepsilon'}(S_T)}
  \lesssim \mathcal I_0
  + M \bignorm{\psi^k_\mp}_{L^2([0,T];H^\varepsilon)}
  \\
  + \sigma(T) \bignorm{\Pi_\pm(D) \bigl(\phi^k
  \beta\psi^k\bigr)}_{X_\pm^{\varepsilon,-1/2+2\varepsilon'}(S_T)},
\end{multline}
\begin{multline}\label{2}
  \bignorm{\phi^{k+1}_\pm}_{\mathcal H^{1/2+\varepsilon,1/2+\varepsilon'}(S_T)}
  \lesssim \mathcal I_0
  + m^2 \bignorm{\phi^k}_{L^2([0,T];H^{-1/2+\varepsilon})}
  \\
  + \sigma(T)
  \bignorm{\biginnerprod{\beta\psi^k}{\psi^k}}_{H^{-1/2+\varepsilon,-1/2+2\varepsilon'}(S_T)},
\end{multline}
so in view of \eqref{0}, the key is to establish the general estimates
\begin{align}
  \label{3}
  \norm{\Pi_\pm(D) \Bigl(\phi
  \beta\Pi_{[\pm]}(D)\psi\Bigr)}_{X_\pm^{\varepsilon,-1/2+2\varepsilon'}}
  &\lesssim \norm{\phi}_{H^{1/2+\varepsilon,1/2+\varepsilon'}}
  \norm{\psi}_{X_{[\pm]}^{\varepsilon,1/2+\varepsilon'}},
  \\
  \label{4}
  \norm{\innerprod{\beta\Pi_{[\pm]}(D)\psi}{\Pi_\pm(D)\psi'}}_{H^{-1/2+\varepsilon,-1/2+2\varepsilon'}}
  &\lesssim \norm{\psi}_{X_{[\pm]}^{\varepsilon,1/2+\varepsilon'}}
  \norm{\psi'}_{X_\pm^{\varepsilon,1/2+\varepsilon'}}
\end{align}
for $\varepsilon'>0$ sufficiently small depending on $\varepsilon$.
Here $\pm$ and $[\pm]$ denote independent signs.
Also, in \eqref{3} and \eqref{4} the norms are
not restricted to $S_T$, but they are applied to extensions of the
iterates $\phi^k,\psi^k$. Taking \emph{infima} over all
extensions, one then obtains, from \eqref{1}--\eqref{4},
$$
  A_{k+1}(T) \lesssim \mathcal I_0 + \sigma(T) P\bigl( A_k(T) \bigr),
$$
where
$$
  A_k(T) = \bignorm{\phi^k}_{\mathcal H^{1/2+\varepsilon,1/2+\varepsilon'}(S_T)}
  + \sum_{\pm} \bignorm{\psi^k_\pm}_{X_\pm^{\varepsilon,1/2+\varepsilon'}(S_T)}
$$
and $P$ is a polynomial such that $P(0) = 0$. Here we used also
the fact that
$$
  \norm{\psi}_{L^2([0,T];H^\varepsilon)}
  \le T^{1/2} \norm{\psi}_{L^\infty([0,T];H^\varepsilon)}
  \lesssim T^{1/2}
  \norm{\psi}_{X_\pm^{\varepsilon,1/2+\varepsilon'}(S_T)},
$$
and similarly for the term
$\norm{\phi^k}_{L^2([0,T];H^{-1/2+\varepsilon})}$ in \eqref{2}.

In view of the multilinearity of the nonlinear terms in the
system, one obtains similar estimates for the difference of
subsequent iterates, and the standard contraction argument then
gives local existence and uniqueness (in the iteration space) for
$0 < T \le 1$ sufficiently small depending on $\mathcal I_0$.

Thus, we have reduced to proving \eqref{3} and \eqref{4}.
But by duality, \eqref{3} is equivalent to
$$  \iint \innerprod{\Pi_\pm(D)\left(\phi\beta\Pi_{[\pm]}(D)\psi\right)}{\psi'} \,
  dt \, dx
  \lesssim  \norm{\phi}_{H^{1/2+\varepsilon,1/2+\varepsilon'}}
  \norm{\psi}_{X_{[\pm]}^{\varepsilon,1/2+\varepsilon'}}
  \norm{\psi'}_{X_\pm^{-\varepsilon,1/2-2\varepsilon'}},
$$
and since
\begin{align*}
  &\iint\innerprod{\Pi_\pm(D)\left(\phi\beta\Pi_{[\pm]}(D)\psi\right)}{\psi'}  \,
  dt \, dx
  \\
  &\qquad\qquad\qquad\qquad = \iint \innerprod{\phi\beta\Pi_{[\pm]}(D)\psi}{\Pi_\pm(D)\psi'} \,
  dt \, dx
  \\
  &\qquad\qquad\qquad\qquad = \iint \phi \innerprod{\beta\Pi_{[\pm]}(D)\psi}{\Pi_\pm(D)\psi'} \,
  dt \, dx
  \\
  &\qquad\qquad\qquad\qquad \le \norm{\phi}_{H^{1/2+\varepsilon,1/2+\varepsilon'}}
  \norm{\innerprod{\beta\Pi_{[\pm]}(D)\psi}{\Pi_\pm(D)\psi'}}_{H^{-1/2-\varepsilon,-1/2-\varepsilon'}},
\end{align*}
we conclude that \eqref{3} reduces to an estimate similar to \eqref{4}:
\begin{equation}\label{6}
  \norm{\innerprod{\beta\Pi_{[\pm]}(D)\psi}{\Pi_\pm(D)\psi'}}_{H^{-1/2-\varepsilon,-1/2-2\varepsilon'}}
  \lesssim \norm{\psi}_{X_{[\pm]}^{\varepsilon,1/2+\varepsilon'}}
  \norm{\psi'}_{X_\pm^{-\varepsilon,1/2-2\varepsilon'}}.
\end{equation}
Note that both \eqref{4} and \eqref{6} are perturbations around
the false estimate \eqref{KeyBilinearProblem}. In fact, using the null structure we
shall reduce \eqref{4} and \eqref{6} to some well-known bilinear
spacetime estimates of Klainerman-Machedon type for products of
free waves. Specifically, we need the following generalization of
the classical $L^4$ estimate \eqref{StrichartzL4} of Strichartz.

\begin{theorem}
\cite{KlainermanMachedon93,
KlainermanMachedon96,FoschiKlainerman00}. Let
$s_1, s_2, s_3 \in \R$. The estimate
$$
  \bignorm{\abs{D}^{-s_3} (uv)}_{L^2(\R^{1+3})} \le
  C_{s_1,s_2,s_3}
  \norm{u_0}_{\dot H^{s_1}} \norm{v_0}_{\dot H^{s_2}}
$$
holds for free waves $u(t) = e^{\pm it\abs{D}}u_0$, $v(t) = e^{\pm
it\abs{D}}v_0$ if and only if
\begin{equation}\label{7}
s_1+s_2+s_3 = 1,
\qquad
s_1,s_2,s_3 < 1,
\qquad
s_1+s_2 > \frac{1}{2}.
\end{equation}
\end{theorem}

From this and Lemma \ref{Lemma2} we obtain
\begin{equation}\label{8}
  H^{s_1,b} \cdot H^{s_2,b} \longrightarrow H^{-s_3,0}
  \quad \text{for $b > 1/2$ and $s_1,s_2,s_3 \ge 0$ satisfying \eqref{7}.}
\end{equation}
Here we use the following notation: If $X,Y,Z$ are normed spaces
of functions, the statement $X \cdot Y \to Z$ means that the
bilinear estimate $\norm{uv}_Z \le C \norm{u}_X \norm{v}_Y$ holds
for some constant $C$.

By interpolation between \eqref{8} and the estimate
\begin{equation}\label{9}
  L^2 \cdot H^{0,b} \longrightarrow H^{-N,0} \qquad \text{for}
  \qquad b > \frac{1}{2}, \, N > \frac{3}{2},
\end{equation}
which by duality is equivalent to $H^{N,0} \cdot H^{0,b} \to L^2$
and therefore follows from H\"older's inequality $L_t^\infty L_x^2
\cdot L_t^2 L_x^\infty \to L^2$ and Sobolev embedding, it is easy
to prove (see the next section) the estimates
\begin{align}
  \label{10}
  H^{0,1/2-\delta} \cdot H^{1/2+\varepsilon,b} &\longrightarrow
  H^{-1/2,0},
  \\
  \label{11}
  H^{0,1/2-\delta} \cdot H^{1/2,b} &\longrightarrow
  H^{-1/2-\varepsilon,0},
  \\
  \label{12}
  H^{1/2-\varepsilon,1/2-\delta} \cdot H^{1/2+\varepsilon,b} &\longrightarrow
  H^{-\varepsilon,0},
  \\
  \label{13}
  H^{1/2-\varepsilon,1/2-\delta} \cdot H^{\varepsilon,b} &\longrightarrow
  H^{-1/2-\varepsilon,0},
  \\
  \label{14}
  H^{-\varepsilon,1/2-\delta} \cdot H^{1/2+\varepsilon,b} &\longrightarrow
  H^{-1/2-\varepsilon,0},
\end{align}
for all $b > 1/2$, $\varepsilon > 0$ and sufficiently small
$\delta > 0$, depending on $\varepsilon$.

We now turn to the proofs of \eqref{4} and \eqref{6}. It suffices to consider the
case where the sign $[\pm]$ is a $+$.
Throughout the rest of this section, we set $b = 1/2+\varepsilon'$; recall that
$\varepsilon'>0$ denotes a sufficiently small number, depending on $\varepsilon$.

\subsection{Proof of \eqref{4}}

Using Lemmas \ref{SymbolLemma} and \ref{Lemma5}, we see that
\eqref{4} reduces to proving
\begin{multline}\label{15}
  \norm{ \frac{1}{\elwt{\xi}^{1/2-\varepsilon}
  \elwt{\abs{\tau}-\abs{\xi}}^{1/2-2\varepsilon'}}
  \iint \theta_{\pm} \bigabs{ \widetilde \psi(\lambda,\eta)}
  \bigabs{ \widetilde{\psi'}(\lambda-\tau,\eta-\xi)}
  \, d\lambda \, d\eta }_{L^2_{\tau,\xi}}
  \\
  \lesssim
  \norm{\psi}_{X_+^{\varepsilon,1/2+\varepsilon'}}
  \norm{\psi'}_{X_\pm^{\varepsilon,1/2+\varepsilon'}},
\end{multline}
where
$$
  \theta_{\pm} \equiv \angle\bigl(\eta,\pm(\eta-\xi) \bigr).
$$
We claim that
$$
  \theta_{+}^2 \approx \frac{\abs{\xi}
  r_{+}}{\abs{\eta}\abs{\eta-\xi}},
  \qquad
  \theta_{-}^2 \approx \frac{(\abs{\eta}+\abs{\eta-\xi})
  r_{-}}{\abs{\eta}\abs{\eta-\xi}},
$$
where
$$
  r_{+} \equiv \abs{\xi} - \bigabs{ \abs{\eta} - \abs{\eta-\xi}
  },
  \qquad
  r_{-} \equiv \abs{\eta} + \abs{\eta-\xi} - \abs{\xi}.
$$
To prove the estimate for $\theta_+^2$, one writes
$$
  \frac{\abs{\xi}
  r_{+}}{\abs{\eta}\abs{\eta-\xi}}
  \approx \frac{\left(\abs{\xi} + \bigabs{ \abs{\eta} - \abs{\eta-\xi}}\right)
  r_{+}}{\abs{\eta}\abs{\eta-\xi}}
  = \frac{\abs{\xi}^2 - \bigabs{ \abs{\eta} - \abs{\eta-\xi}}^2}{\abs{\eta}\abs{\eta-\xi}}
  = 2 \left( 1 - \cos \theta_+ \right),
$$
and the estimate for $\theta_-^2$ is proved in a similar way.

Assuming as we may that $\widetilde \psi, \widetilde{\psi'}
\ge 0$, and suppressing the arguments of these functions to keep
the notation manageable, we thus reduce \eqref{15} to
\begin{multline}
\label{16}
  \norm{ \frac{\abs{\xi}^\varepsilon}{\elwt{\abs{\tau}-\abs{\xi}}^{1/2-2\varepsilon'}}
  \iint \frac{r_{+}^{1/2}}{\abs{\eta}^{1/2}\abs{\eta-\xi}^{1/2}}
  \widetilde \psi
  \widetilde{\psi'}
  \, d\lambda \, d\eta }_{L^2_{\tau,\xi}}
  \\
  \lesssim
  \norm{\psi}_{X_+^{\varepsilon,1/2+\varepsilon'}}
  \norm{\psi'}_{X_+^{\varepsilon,1/2+\varepsilon'}}
\end{multline}
and
\begin{multline}
  \label{17}
  \norm{ \frac{1}{\elwt{\xi}^{1/2-\varepsilon}
  \elwt{\abs{\tau}-\abs{\xi}}^{1/2-2\varepsilon'}}
  \iint \frac{r_{-}^{1/2}}{\min\bigl(\abs{\eta},\abs{\eta-\xi}\bigr)^{1/2}}
  \widetilde \psi
  \widetilde{\psi'}
  \, d\lambda \, d\eta }_{L^2_{\tau,\xi}}
  \\
  \lesssim
  \norm{\psi}_{X_+^{\varepsilon,1/2+\varepsilon'}}
  \norm{\psi'}_{X_-^{\varepsilon,1/2+\varepsilon'}}.
\end{multline}

Now we apply the following:

\begin{lemma}\label{Lemma6}
$r_{\pm} \lesssim \bigabs{\abs{\tau}-\abs{\xi}} +
\bigabs{\lambda+\abs{\eta}} +
\bigabs{\lambda-\tau\pm\abs{\eta-\xi}}$
\end{lemma}

\begin{proof}
If $\tau \ge 0$, we estimate
$$
  r_{+} \le \abs{\xi} + \abs{\eta} - \abs{\eta-\xi} = \abs{\xi} -
  \abs{\tau} + \lambda + \abs{\eta} + \tau - \lambda -
  \abs{\eta-\xi},
$$
while if $\tau < 0$ we use
$$
  r_{+} \le \abs{\xi} - \abs{\eta} + \abs{\eta-\xi} = \abs{\xi} -
  \abs{\tau} - \lambda - \abs{\eta} - \tau + \lambda +
  \abs{\eta-\xi}.
$$
	To handle $r_{-}$ we write
$$
  r_{-} = \lambda + \abs{\eta} + \tau - \lambda + \abs{\eta-\xi} - \tau
  - \abs{\xi}.
$$
If $\tau < 0$, this equals $\lambda + \abs{\eta} + \tau - \lambda
+ \abs{\eta-\xi} + \abs{\tau} - \abs{\xi}$, while if $\tau \ge 0$
it is $\le \lambda + \abs{\eta} + \tau - \lambda +
\abs{\eta-\xi}$.
\end{proof}

We also need
\begin{equation}\label{18}
  r_{\pm} \le 2 \min\bigl(\abs{\eta},\abs{\eta-\xi}\bigr),
\end{equation}
which follows from the triangle inequality. Combinining this with
Lemma \ref{Lemma6} we get
$$
  r_{\pm}^{1/2} \lesssim \bigabs{\abs{\tau}-\abs{\xi}}^{1/2-2\varepsilon'}
  \min\bigl(\abs{\eta},\abs{\eta-\xi}\bigr)^{2\varepsilon'}
  + \bigabs{\lambda+\abs{\eta}}^{1/2}
  + \bigabs{\lambda-\tau\pm\abs{\eta-\xi}}^{1/2}.
$$
Moreover, by symmetry we may assume $\abs{\eta} \ge
\abs{\eta-\xi}$ in \eqref{16} and \eqref{17}. Hence \eqref{16}
reduces to proving
$$
  I_j^{+} \lesssim
  \norm{\psi}_{X_+^{\varepsilon,1/2+\varepsilon'}}
  \norm{\psi'}_{X_+^{\varepsilon,1/2+\varepsilon'}}
$$
for $j = 1,2,3$, where
\begin{align*}
  I_1^{+} &= \norm{
  \iint \frac{\widetilde \psi
  \widetilde{\psi'}}{\abs{\eta}^{1/2-\varepsilon}\abs{\eta-\xi}^{1/2-2\varepsilon'}}
  \, d\lambda \, d\eta }_{L^2_{\tau,\xi}},
  \\
  I_2^{+} &= \norm{
  \frac{1}{\elwt{\abs{\tau}-\abs{\xi}}^{1/2-2\varepsilon'}}
  \iint \frac{\bigabs{\lambda+\abs{\eta}}^{1/2} \widetilde \psi
  \widetilde{\psi'}}{\abs{\eta}^{1/2-\varepsilon}\abs{\eta-\xi}^{1/2}}
  \, d\lambda \, d\eta }_{L^2_{\tau,\xi}},
  \\
  I_3^{+} &= \norm{
  \frac{1}{\elwt{\abs{\tau}-\abs{\xi}}^{1/2-2\varepsilon'}}
  \iint \frac{\widetilde \psi\,
  \bigabs{\lambda-\tau+\abs{\eta-\xi}}^{1/2}\widetilde{\psi'}}
  {\abs{\eta}^{1/2-\varepsilon}\abs{\eta-\xi}^{1/2}}
  \, d\lambda \, d\eta }_{L^2_{\tau,\xi}},
\end{align*}
and \eqref{17} reduces to proving
$$
  I_j^{-} \lesssim
  \norm{\psi}_{X_+^{\varepsilon,1/2+\varepsilon'}}
  \norm{\psi'}_{X_-^{\varepsilon,1/2+\varepsilon'}}
$$
for $j = 1,2,3$, where
\begin{align*}
  I_1^{-} &= \norm{ \frac{1}{\elwt{\xi}^{1/2-\varepsilon}}
  \iint \frac{\widetilde \psi
  \widetilde{\psi'}}{\abs{\eta-\xi}^{1/2-2\varepsilon'}}
  \, d\lambda \, d\eta }_{L^2_{\tau,\xi}},
  \\
  I_2^{-} &= \norm{ \frac{1}{\elwt{\xi}^{1/2-\varepsilon}
  \elwt{\abs{\tau}-\abs{\xi}}^{1/2-2\varepsilon'}}
  \iint \frac{\bigabs{\lambda+\abs{\eta}}^{1/2} \widetilde \psi
  \widetilde{\psi'}}{\abs{\eta-\xi}^{1/2}}
  \, d\lambda \, d\eta }_{L^2_{\tau,\xi}},
  \\
  I_3^{-} &= \norm{ \frac{1}{\elwt{\xi}^{1/2-\varepsilon}
  \elwt{\abs{\tau}-\abs{\xi}}^{1/2-2\varepsilon'}}
  \iint \frac{\widetilde \psi\,
  \bigabs{\lambda-\tau-\abs{\eta-\xi}}^{1/2}\widetilde{\psi'}}{\abs{\eta-\xi}^{1/2}}
  \, d\lambda \, d\eta }_{L^2_{\tau,\xi}}.
\end{align*}

The estimates for $I_j^{+}$, $j = 1,2,3$, reduce to,
respectively,
\begin{align}
  \label{19}
  X_+^{1/2,b} \cdot X_+^{1/2+\varepsilon-2\varepsilon',b}
  &\longrightarrow L^2,
  \\
  \label{20}
  X_+^{1/2,0} \cdot X_+^{1/2+\varepsilon,b}
  &\longrightarrow H^{0,-1/2+2\varepsilon'},
  \\
  \label{21}
  X_+^{1/2,b} \cdot X_+^{1/2+\varepsilon,0}
  &\longrightarrow H^{0,-1/2+2\varepsilon'},
\end{align}
and the estimates for $I_j^{-}$, $j = 1,2,3$, reduce to
\begin{align}
  \label{22}
  X_+^{\varepsilon,b} \cdot X_-^{1/2+\varepsilon-2\varepsilon',b}
  &\longrightarrow H^{-1/2+\varepsilon,0},
  \\
  \label{23}
  X_+^{\varepsilon,0} \cdot X_-^{1/2+\varepsilon,b}
  &\longrightarrow H^{-1/2+\varepsilon,-1/2+2\varepsilon'},
  \\
  \label{24}
  X_+^{\varepsilon,b} \cdot X_-^{1/2+\varepsilon,0}
  &\longrightarrow H^{-1/2+\varepsilon,-1/2+2\varepsilon'}.
\end{align}
Recall that we assume $b = 1/2+\varepsilon'$ throughout this section.

Since in general $\norm{u}_{H^{s,b}} \le \norm{u}_{X_\pm^{s,b}}$ for $b \ge 0$,
we may in fact replace all the $X^{s,b}_\pm$ in the left hand
sides by $H^{s,b}$. (The information encoded in the signs has already been made use of
through the null structure.) Using also duality, we thus reduce to
\begin{align*}
  H^{1/2,b} \cdot H^{1/2+\varepsilon-2\varepsilon',b}
  &\longrightarrow L^2,
  \\
  H^{0,1/2-2\varepsilon'} \cdot H^{1/2+\varepsilon,b}
  &\longrightarrow H{-1/2,0},
  \\
  H^{1/2,b} \cdot H^{0,1/2-2\varepsilon'}
  &\longrightarrow H^{-1/2-\varepsilon,0},
\end{align*}
and
\begin{align*}
  H^{\varepsilon,b} \cdot H^{1/2+\varepsilon-2\varepsilon',b}
  &\longrightarrow H^{-1/2+\varepsilon,0},
  \\
  H^{1/2-\varepsilon,1/2-2\varepsilon'}  \cdot H^{1/2+\varepsilon,b}
  &\longrightarrow H^{-\varepsilon,0},
  \\
  H^{\varepsilon,b} \cdot H^{1/2-\varepsilon,1/2-2\varepsilon'}
  &\longrightarrow H^{-1/2-\varepsilon,0}.
\end{align*}
All these estimates are true for $\varepsilon'>0$ sufficiently
small, by \eqref{8} and \eqref{10}--\eqref{14}.

\subsection{Proof of \eqref{6}}

Proceeding as in the proof of \eqref{4}, we reduce to
\begin{multline}
\label{25}
  \norm{ \frac{\abs{\xi}^\varepsilon}{\elwt{\abs{\tau}-\abs{\xi}}^{1/2+\varepsilon'}}
  \iint \frac{r_{+}^{1/2}}{\abs{\eta}^{1/2}\abs{\eta-\xi}^{1/2}}
  \widetilde \psi
  \widetilde{\psi'}
  \, d\lambda \, d\eta }_{L^2_{\tau,\xi}}
  \\
  \lesssim
  \norm{\psi}_{X_+^{\varepsilon,1/2+\varepsilon'}}
  \norm{\psi'}_{X_+^{-\varepsilon,1/2-2\varepsilon'}}
\end{multline}
and
\begin{multline}
  \label{26}
  \norm{ \frac{1}{\elwt{\xi}^{1/2+\varepsilon}
  \elwt{\abs{\tau}-\abs{\xi}}^{1/2+\varepsilon'}}
  \iint \frac{r_{-}^{1/2}}{\min\bigl(\abs{\eta},\abs{\eta-\xi}\bigr)^{1/2}}
  \widetilde \psi
  \widetilde{\psi'}
  \, d\lambda \, d\eta }_{L^2_{\tau,\xi}}
  \\
  \lesssim
  \norm{\psi}_{X_+^{\varepsilon,1/2+\varepsilon'}}
  \norm{\psi'}_{X_-^{-\varepsilon,1/2-2\varepsilon'}}.
\end{multline}
By Lemma \ref{Lemma6} and \eqref{18},
$$
  r_{\pm}^{1/2} \lesssim \bigabs{\abs{\tau}-\abs{\xi}}^{1/2}
  + \bigabs{\lambda+\abs{\eta}}^{1/2}
  + \bigabs{\lambda-\tau\pm\abs{\eta-\xi}}^{1/2-2\varepsilon'}
  \min\bigl(\abs{\eta},\abs{\eta-\xi}\bigr)^{2\varepsilon'},
$$
hence \eqref{25} reduces to (recall that $b = 1/2+\varepsilon'$)
\begin{align*}
  H^{1/2+\varepsilon,b} \cdot H^{1/2-\varepsilon,1/2-\varepsilon'}
  &\longrightarrow H^{-\varepsilon,0},
  \\
  H^{1/2+\varepsilon,0} \cdot H^{1/2-\varepsilon,1/2-\varepsilon'}
  &\longrightarrow H^{-\varepsilon,-b},
  \\
  H^{1/2+\varepsilon-2\varepsilon',b} \cdot H^{1/2-\varepsilon-2\varepsilon',0}
  &\longrightarrow H^{-\varepsilon,-b},
\end{align*}
which follow from, respectively, \eqref{12}, \eqref{13} (via
duality) and \eqref{8} (also via duality). Now consider
\eqref{26}. Assuming first $\abs{\eta} \le \abs{\eta-\xi}$ we
reduce to
\begin{align*}
  H^{1/2+\varepsilon,b} \cdot H^{-\varepsilon,1/2-2\varepsilon'}
  &\longrightarrow H^{-1/2-\varepsilon,0},
  \\
  H^{1/2+\varepsilon,0} \cdot H^{-\varepsilon,1/2-2\varepsilon'}
  &\longrightarrow H^{-1/2-\varepsilon,-b},
  \\
  H^{1/2+\varepsilon-2\varepsilon',b} \cdot H^{-\varepsilon,0}
  &\longrightarrow H^{-1/2-\varepsilon,-b},
\end{align*}
while in the case $\abs{\eta} \ge \abs{\eta-\xi}$ we get
\begin{align*}
  H^{\varepsilon,b} \cdot H^{1/2-\varepsilon,1/2-2\varepsilon'}
  &\longrightarrow H^{-1/2-\varepsilon,0},
  \\
  H^{\varepsilon,0} \cdot H^{1/2+-\varepsilon,1/2-2\varepsilon'}
  &\longrightarrow H^{-1/2-\varepsilon,-b},
  \\
  H^{\varepsilon,b} \cdot H^{1/2-\varepsilon-2\varepsilon',0}
  &\longrightarrow H^{-1/2-\varepsilon,-b}.
\end{align*}
All these reduce (possibly via duality) to \eqref{8} or
\eqref{10}--\eqref{14}.

\section{Proof of \eqref{10}--\eqref{14}}

All these follow by interpolation between \eqref{9} and various
special cases of \eqref{8}. Fix $\varepsilon > 0$, $b > 1/2$, $N
> 3/2$. The number $\delta > 0$ will be chosen sufficiently small,
depending on $\varepsilon$.

For \eqref{10} we interpolate between
\begin{align*}
  H^{0,1/2+\delta} \cdot H^{1/2+\varepsilon,b} &\longrightarrow
  H^{-1/2+\varepsilon,0},
  \\
  L^2 \cdot H^{1/2+\varepsilon,b} &\longrightarrow
  H^{-N,0}.
\end{align*}
This gives
$$
  H^{0,(1-\theta)(1/2+\delta)} \cdot H^{1/2+\varepsilon,b} \longrightarrow
  H^{(1-\theta)(-1/2+\varepsilon) - \theta N,0}
$$
for $0 \le \theta \le 1$. First choose $\theta > 0$ so small that
$(1-\theta)(-1/2+\varepsilon) - \theta N \ge - 1/2$. Then choose
$\delta > 0$ so small that $(1-\theta)(1/2+\delta) \le 1/2 -
\delta$.

For \eqref{11} we interpolate between
\begin{align*}
  H^{0,1/2-\delta} \cdot H^{1/2+\varepsilon',b} &\longrightarrow
  H^{-1/2,0},
  \\
  H^{0,1/2-\delta} \cdot H^{0,b} &\longrightarrow
  H^{-N,0},
\end{align*}
the first of which is just \eqref{10}, for $\varepsilon',\delta >
0$ to be chosen sufficiently small, depending on $\varepsilon$.
This gives
$$
  H^{0,1/2-\delta} \cdot H^{(1-\theta)(1/2+\varepsilon'),b} \longrightarrow
  H^{(1-\theta)(-1/2) - \theta N,0}.
$$
First choose $\theta > 0$ so small that $(1-\theta)(-1/2) - \theta
N \ge -1/2 - \varepsilon$. Then choose $\varepsilon' > 0$ so small
that $(1-\theta)(1/2+\varepsilon') \le 1/2$.

For \eqref{12} interpolate between
\begin{align*}
  H^{1/2-\varepsilon,1/2+\delta} \cdot H^{1/2+\varepsilon,b} &\longrightarrow
  L^2,
  \\
  H^{1/2-\varepsilon,0} \cdot H^{1/2+\varepsilon,b} &\longrightarrow
  H^{-N,0}.
\end{align*}
Thus,
$$
  H^{1/2-\varepsilon,(1-\theta)(1/2+\delta)} \cdot H^{1/2+\varepsilon,b} \longrightarrow
  H^{- \theta N,0}
$$
for $0 \le \theta \le 1$. First choose $\theta > 0$ so small that
$\theta N \le \varepsilon$. Then choose $\delta > 0$ so small that
$(1-\theta)(1/2+\delta) \le 1/2 - \delta$.

To prove \eqref{13} we interpolate
\begin{align*}
  H^{1/2-\varepsilon,1/2+\delta} \cdot H^{\varepsilon,b} &\longrightarrow
  H^{-1/2,0},
  \\
  H^{1/2-\varepsilon,0} \cdot H^{\varepsilon,b} &\longrightarrow
  H^{-N,0},
\end{align*}
yielding
$$
  H^{1/2-\varepsilon,(1-\theta)(1/2+\delta)} \cdot H^{\varepsilon,b} \longrightarrow
  H^{(1-\theta)(-1/2) - \theta N,0}.
$$
First choose $\theta > 0$ so small that $(1-\theta)(1/2) + \theta
N \le 1/2+ \varepsilon$. Then choose $\delta > 0$ so small that
$(1-\theta)(1/2+\delta) \le 1/2 - \delta$.

Finally, \eqref{14} reduces to \eqref{11} or \eqref{12}, depending
on whether the frequency interactions are of type high-high or
high-low/low-high in the product on the left hand side.

\end{document}